\documentstyle{amsppt} \tolerance 3000 \pagewidth{5.5in}
\vsize7.0in \magnification=\magstep1 \widestnumber
\key{AAAAAAAAAAAAAAAAA} \NoRunningHeads \topmatter
\title Distinct distances on a sphere
\endtitle
\author Alex Iosevich and Mischa Rudnev
\endauthor
\date August 18, 2004
\enddate
\address Department of Mathematics, University of Missouri-Columbia,
Columbia MO 65211 USA \endaddress \email iosevich \@
math.missouri.edu http://www.math.missouri.edu/ $\sim$ iosevich
\endemail
\address Department of Mathematics, University of Bristol,
Bristol BS8 1TW UK \endaddress \email m.rudnev \@ bris.ac.uk
\endemail
\abstract We prove that a set of $N$ points on a two dimensional
sphere satisfying a discrete energy condition determines at least
a constant times $N$ distinct distances. Homogeneous sets in the
sense of Solymosi and Vu easily satisfy this condition, as do
other sets that in the sense that will be made precise below
respect the curvature properties of the sphere. \endabstract
\thanks The work was partly supported by a grant from the National
Science Foundation NSF02-45369
\endthanks
\endtopmatter
\document

The classical Erd\"os distance conjecture (EDC) says that a planar
point set of cardinality $N$ determines at least a constant times
$\frac{N}{\sqrt{\log(N)}}$ Euclidean distances. Taking
$A={[0,\sqrt{N}]}^2 \cap {\Bbb Z}^2$ shows that such an estimate
would be best possible (\cite{Erd46}). More precisely, let $A
\subset {\Bbb R}^2$ with $\# A=N$. Let
$$\Delta(A)=\{\|x-y\|: a,b \in A\}, \tag0.1$$
where $\|z\|=\sqrt{z_1^2+z_2^2}$. See, for example, \cite{PA96}
and \cite{PS04} for a description of this beautiful problem and
connections with other problems in geometric combinatorics.

In spite of many efforts over the past sixty years, the problem
remains unsolved. The best known result to date, due to Katz and
Tardos (\cite{KT04}), partly based on an ingenuous combinatorial
technique by Solymosi and T\'oth (\cite{ST01}), gives
$$ \# \Delta(A) \gtrsim N^{\approx .86}. \tag0.2$$

Here and throughout the paper, $a \lesssim b$ ($a \gtrsim b$)
means that there exists a universal constant $C$ such that $a \leq
Cb$ ($a\geq Cb$) and $a \approx b$ means that $a\lesssim b$ and
$b\lesssim a.$

Many of the aforementioned papers on the EDC observe that their
method extends to distinct distances on spheres. This is partly
due to the fact that many of the principles of planar geometry
such as the existence of easily constructed geodesics and
bisectors carries over to the spherical setting. In this paper we
shall see that the curvature properties of the unit sphere allow
for a very efficient accounting procedure for the distance set and
lead to a complete resolution of the distance conjecture in this
context, under an additional assumption on the set $A$, that in
effect enables one to take advantage of curvature. Our main result
is the following.

\proclaim{Theorem 0.1} Let $A \subset S^2$ be a finite set of
cardinality $N$. Let $\Delta(A)$ be defined as above. Let
$\theta(a,b)$ denote the spherical, or angular, distance between
$a$ and $b$, $a,b \in S^2$. Then
$$ \# \Delta(A) \gtrsim \frac{N}{I_1(A)}, \tag0.3$$ where
$$ I_{\beta}(N)=\frac{1}{N^2} \sum_{a\not=b}
\frac{1}{\theta^{\beta}(a,b)}. \tag0.4$$

\endproclaim

The applicability of Theorem 0.1 depends on the extent to which we
can bound $I_1(N)$ from above. For example, the quantity $I_1(N)$
is unbounded, as $N \to \infty$, if $A$ is a maximal
$\frac{1}{N}$-separated subset of $S^2$ contained in a curved
$\frac{1}{\sqrt{N}}$ by $\frac{1}{\sqrt{N}}$ rectangle. On the
other hand, $I_1(N)$ is bounded if $A$ is $\frac{1}{\sqrt{N}}$
separated. Observe that one can partition the sphere into $N$
$\frac{1}{\sqrt{N}}$ by $\frac{1}{\sqrt{N}}$ curved rectangles
with the property that any point of $S^2$ is covered by at most
three of these rectangles. Thus this special case of boundedness
of $I_1(N)$ can be viewed as a spherical analog of the homogeneity
condition on the point set employed by Solymosi and Vu
(\cite{SV03}). Also observe that if $A$ consists of $N$ equally
spaced points on a great circle in $S^2$, $I_1(N)$ grows
logarithmically in $N$.

We conclude that for a large class of point sets on $S^2$, the EDC
holds. Moreover, we conjecture that on $S^2$, the Erd\"os Distance
Conjecture should hold without the logarithmic factor as the
lattice example from ${\Bbb R}^2$ has no apparent analog in this
setting.

Observe that $I_{\beta}(A)$ is a discrete analog of the energy
integral used to measure Hausdorff dimension of sets in a
continuous setting. This connection is explored at the end of the
paper.

The method of proof of Theorem 0.1 easily generalizes to higher
dimensional spheres. In the process, interesting issues involving
the geometric properties of the function $I_1(A)$ arise. We shall
investigate this matter systematically in a subsequent paper.

\vskip.125in

\head Proof of Theorem 0.1 \endhead

\vskip.125in

Let $a,b \in S^2$. Then ${\|a-b\|}^2={\|a\|}^2+{\|b\|}^2-2a \cdot
b=2-2a \cdot b$. It follows that instead of counting Euclidean
distances on $S^2$, it suffices to count dot products $a \cdot b$.
Let $h_{a, \delta}$ be a smooth cutoff function on $S^2$,
identically equal to $1$ in the $\delta$-neighborhood of $a \in
S^2$, and vanishing outside the $2 \delta$-neighborhood of $a$.
For convenience we construct this function in such a way that any
$h_{a, \delta}$ can be obtained from a $h_{b, \delta}$ by a
rotation that takes $a$ to $b$. The right choice for $\delta$ will
turn out to be $\frac{1}{N}$, but we will keep this parameter
flexible for a while for the sake of clarity.

\subhead Construction of measures approximating $A$ and
$\Delta(A)$ \endsubhead Let $\omega\in S^2$ and
$$ \mu(\omega)=c_1\frac{1}{N \delta^2} \sum_{a \in A}
h_{a, \delta}(\omega), \tag1.1$$ a measure approximating $A$.
Clearly the constant $c_1$ can be chosen such that $\mu$ is a
probability measure, so in the sequel assume $c_1=1$. We now
construct a measure approximating $\Delta(A)$. Let $A_{\delta}$
denote the support of $\mu$, the set where $\mu$ is non-zero. Let
$f$ be a function on $\Delta(A_{\delta})$. Then the distance
measure $\nu$ is defined via the following identity:
$$ \int f(t) d\nu(t)=\int \int f(\|x-y\|) d\mu(x)d\mu(y). \tag1.2$$
In case of $S^2$, see above, when one looks at angular distances,
instead of (1.2) one can use the following definition for the
distance measure $\nu$:
$$ \int f(t) d\nu(t)=\int \int f(x\cdot y) d\mu(x)d\mu(y). \tag1.3$$
Our plan is to show that $\nu$ defined via (1.3) has an $L^2$
density and then convert this statement into a lower bound for the
number of distances. On the combinatorial level, what we are doing
can be described as follows. Let $m(t)=\# \{(x,y) \in A \times A:
\|x-y\|=t\}$, the incidence function. Estimating the $L^2$ norm of
$\nu$ with an appropriate $\delta$, turns out to be equivalent to
estimating $\sum_{t \in \Delta(A)} m^2(t)$ from above. The latter
is easily converted into a lower bound for the number of distances
using the Cauchy-Schwartz inequality and the fact that $\sum_{t
\in \Delta(A)} m(t) \approx N^2$.

\subhead Estimation of the $L^2$ norm of $\nu$ \endsubhead Taking
$f(t)=e^{2 \pi i \lambda t}$ in $(1.3)$ we get
$$ \widehat{\nu}(\lambda)=\frac{1}{N^2 \delta^4} \sum_{a,b \in A}
\int_{S^2} \widehat{h_{a,\delta}}(\lambda \omega)h_{b,
\delta}(\omega)d\omega. \tag1.4$$

By the standard method of stationary phase (see e.g. \cite{W03}),
$$ |\widehat{h_{a,\delta}}(\lambda \omega)| \lesssim \min
\{\delta^2, \lambda^{-1}\}, \tag1.5$$ if $\omega \in
U_{a,2\delta}$, where $U_{a, 2\delta}$ is the
$2\delta$-neighborhood of $a$, and
$$ |\widehat{h_{a,\delta}}(\lambda \omega)| \lesssim
\min \{\delta^2, \delta^2 {[\lambda \theta(\omega, U_{a,
\delta})]}^{-M}\}, \tag1.6$$ if $\omega \notin U_{a, 2 \delta}$,
where $\theta(\omega, U_{a, \delta})$ is the angular distance from
$\omega$ to $U_{a, \delta}$. The exponent $M$ in $(1.6)$ is an
arbitrary positive integer,\footnote{The formulas $(1.5)$ and
$(1.6)$ are obtained as follows. If the direction $\omega$ on the
``Fourier side'' coincides with one of the normal directions to
the spherical cap in question, one gets the standard decay for the
Fourier transform of Lebesgue measure on the unit sphere in
$(1.5)$. Otherwise, one can integrate by parts $M$ times and using
the trivial area bound $\delta^2$ for the oscillatory integral in
the final step.} but we shall confine ourselves to the case $M=2$.

Plugging $(1.5)$ and $(1.6)$ into $(1.4)$ yields
$$ |\widehat{\nu}(\lambda)| \lesssim \frac{1}{N}
\chi_{[-\delta^{-2}, \delta^{-2}]}(\lambda)+\frac{1}{N \lambda
\delta^2} \chi_{{\Bbb R} \backslash [-\delta^{-2},
\delta^{-2}]}(\lambda)$$ $$+\frac{1}{N^2} \sum_{a \not=b}
\chi_{[-\theta^{-1}(a,b),
\theta^{-1}(a,b)]}(\lambda)+\frac{1}{N^2} \sum_{a \not=b}
\chi_{{\Bbb R} \backslash [-\theta^{-1}(a,b),
\theta^{-1}(a,b)]}(\lambda){[\lambda \theta(a,b)]}^{-2}$$
$$=I+II+III+IV, \tag1.7$$ where $\chi_J$ is the characteristic
function of a set $J$. Observe that the first line in $(1.7)$ in
essence corresponds to setting $a=b$ in $(1.4)$, when the estimate
$(1.5)$ comes into play, while the second line in $(1.7)$ is the
case $a\neq b$ in $(1.4)$, which uses the estimate $(1.6)$.

A straightforward calculation shows that
$$ \int I^2 d\lambda+ \int {II}^2 d\lambda \lesssim 1, \tag1.8$$
provided that $$\delta \gtrsim \frac{1}{N}.\tag1.9$$

We now estimate $III$ and $IV$. We have
$$ \int {III}^2 d\lambda=\frac{1}{N^4} \sum_{a \not=a'; b
\not=b'} \chi_{[-\theta^{-1}(a,b), \theta^{-1}(a,b)]}(\lambda)
\chi_{[-\theta^{-1}(a',b'),
\theta^{-1}(a',b')]}(\lambda)d\lambda$$ $$ \lesssim \frac{1}{N^2}
\sum_{a \not=b} \theta^{-1}(a,b)=I_1(N). \tag1.10$$ by assumption.
The estimate on $IV$ is clearly the same.

\subhead The lower bound on the cardinality of the distance set
\endsubhead By construction we have that
$$ \int d\nu(t)=1. \tag1.11$$

By Cauchy-Schwartz it follows that
$$ 1 \leq |supp(\nu)| \cdot \int \nu^2(t)dt=
|supp(\nu)| \cdot \int {|\widehat{\nu}(\lambda)|}^2 d\lambda,
\tag1.12$$ where the second equality follows by Plancherel. Note
that $|supp(\nu)|$ denotes the Lebesgue measure of the support of
$\nu$.

By construction, $supp(\nu)$ consists of intervals of length
$\approx \delta$, so
$$ \# \Delta(A) \gtrsim \frac{|supp(\nu)|}{\delta} \gtrsim
\frac{N}{I_1(N)}, \tag1.13$$ since we have chosen $\delta \approx
\frac{1}{N}$.

This completes the proof of Theorem 0.1. We conclude the argument
by pointing out that the quantity
$$I_1(N) \approx \int \frac{d\nu(t)}{t}=\int \int
{|x-y|}^{-1} d\mu(x)d\mu(y), \tag1.14$$ by definition of $\nu$.
The last integral is always finite if $\mu$ uniformly approximates
a measure on a set of Hausdorff dimension greater than one. This
need not be the case, however. 

\newpage

\head References \endhead

\vskip.125in

\ref \key Erd46 \by P. Erd\"os \paper On sets of distances of $n$
points \jour Amer. Math. Monthly \vol 53 \yr 1946 \pages 248-250
\endref

\ref \key KT04 \by N. Katz and G. Tardos \paper A new entropy
inequality for the Erd\"os distance problem \jour Towards a Theory
of Geometric Graphs. (ed.J Pach) Contemporary Mathematics \vol 342
\yr 2004 \endref

\ref \key PA95 \by J. Pach and P. Agarwal \paper Combinatorial
geometry \jour Wiley-Interscience Series in Discrete Mathematics
and Optimization. A Wiley-Interscience Publication. John Wiley and
Sons, Inc., New York \yr 1995 \endref

\ref \key PS04 \by J. Pach and M. Sharir \paper Geometric
incidences \jour Towards a Theory of Geometric Graphs. (ed.J Pach)
Contemporary Mathematics \vol 342 \yr 2004 \endref

\ref \key ST01 \by J. Solymosi and C. T\'oth \paper Distinct
distances in the plane \jour Discr. Comp. Jour. (Misha Sharir
birthday issue) \vol 25 \yr 2001 \pages 629-634 \endref

\ref \key SV03 \by J. Solymosi and V. Vu \paper Distinct distances
in homogeneous sets \jour Symposium on Computational Geometry \yr
2003 \pages 104-105 \endref

\ref \key W03 \by T. Wolff \paper Lectures on harmonic analysis
\jour University Lecture Series, American Mathematical Society,
Providence, RI \vol 29 \yr 2003 \endref

\enddocument